\documentclass{article}%
\usepackage{amsmath}
\usepackage{amsfonts}
\usepackage{amssymb}
\usepackage{graphicx}%
\newtheorem{theorem}{Theorem}[section]

\newtheorem{definition}[theorem]{Definition}
\newtheorem{example}[theorem]{Example}

\newtheorem{lemma}[theorem]{Lemma}

\newtheorem{remark}[theorem]{Remark}

\newenvironment{proof}[1][Proof]{\noindent\textbf{#1.} }{\ \rule{0.5em}{0.5em}}
\begin{document}

\title{Invariance of Arens module actions on the norm dual of a Riesz-normed $f$-algebra}
\author{Karim Boulabiar\\{\small Laboratoire de Recherche LATAO}\\{\small D\'{e}partement de Math\'{e}matiques, Facult\'{e} des Sciences de
Tunis}\\{\small Universit\'{e} de Tunis El Manar, 2092, El Manar, Tunisia}\\{\small Mathematical Institute, Leiden University}\\{\small P.O. Box 9512, 2300 RA Leiden, The Netherlands}\\{\small Email: \texttt{karim.boulabiar@fst.utm.tn}}
\and Rania Souissi\\{\small Laboratoire de Recherche LATAO}\\{\small D\'{e}partement de Math\'{e}matiques, Facult\'{e} des Sciences de
Tunis}\\{\small Universit\'{e} de Tunis El Manar, 2092, El Manar, Tunisia}\\{\small Email: \texttt{rania.souissi@ipest.utm.tn}}}
\date{}
\maketitle

\begin{abstract}
Let $A$ be an $f$-algebra and let $V$ be an order ideal of its order dual
$A^{\sim}$ that separates the points of $A$. In a comment accompanying a
theorem of Huijsmans and de~Pagter (J.~Funct.\ Anal.\ 59 (1984), 41--64) it is
stated that the proof of that theorem can easily be adapted to show that the
order continuous dual $V_{n}^{\sim}$ of $V$ is an $f$-algebra under the Arens
multiplication. We show that this statement is false. The Arens multiplication
on $V_{n}^{\sim}$ requires two invariance conditions on $V$, with respect to
the module actions of $A$ and of $V_{n}^{\sim}$, which do not follow from the
stated hypotheses and can fail already in the most natural case, where $A$
carries a Riesz norm and $V=A^{\ast}$ is the norm dual of $A$. We prove that
the first condition holds if and only if $A$ is bounded, and that the second
holds whenever the Riesz norm is submultiplicative (which forces $A$ to be
bounded); in that case $(A^{\ast})_{n}^{\sim}$ is indeed an $f$-algebra. We
then construct a bounded $f$-algebra of finitely supported real sequences,
equipped with a non-submultiplicative Riesz norm, for which the second
condition fails, so that the Arens multiplication on $(A^{\ast})_{n}^{\sim}$
is not even defined. The comment has also been used in the literature to
derive a further result, which is false as well. This clarifies the exact
scope of the comment.

\end{abstract}


\bigskip

\section{Introduction}

The order dual of a lattice-ordered algebra carries a natural module structure
induced by the multiplication of the algebra. This structure is the starting
point for extending the multiplication to order duals and, following Arens
\cite{A}, for defining a product on order biduals. A basic issue in this
construction is whether the dual spaces involved are invariant under the
induced module actions. Such invariance is needed for the Arens multiplication
to be well defined and, ultimately, for obtaining a lattice-ordered algebra
structure on the order continuous dual.

In their influential and pioneering paper \cite{HdP}, Huijsmans and de Pagter
proved that if $A$ is an $f$-algebra, then the order continuous dual $\left(
A^{\sim}\right)  _{n}^{\sim}$ of the order dual $A^{\sim}$ is an $f$-algebra
under the Arens multiplication (Theorem 4.4 of \cite{HdP}). In a comment
accompanying this theorem, they state that its proof can easily be adapted
when $A^{\sim}$ is replaced by an order ideal $V$ of $A^{\sim}$ that separates
the points of $A$. The proof of Theorem 4.4 is correct, but the comment is
false. The Arens multiplication on $V_{n}^{\sim}$ requires invariance
properties of $V$ that do not follow from the assumption that $V$ is an order
ideal separating the points of $A$, and these properties can fail already in
the most natural case, where $A$ carries a Riesz norm and $V$ is its norm dual
$A^{\ast}$, which is known to be an order ideal of $A^{\sim}$ that separates
the points of $A$: the first property fails as soon as $A$ is not bounded, and
the second one can fail even when $A$ is bounded. The purpose of this note is
to examine these properties and to determine precisely when they hold. The
comment has also been used in the literature for $V=A^{\ast}$ to derive a
further result, which is false as well, as shown in \cite{BS}. This gives an
additional reason to settle exactly when the Arens multiplication on $\left(
A^{\ast}\right)  _{n}^{\sim}$ is available.

To make this precise, let $A$ be an $f$-algebra and let $V$ be a Riesz
subspace of $A^{\sim}$. For $a\in A$ and $f\in V$, define a linear functional
$f\cdot a$ on $A$ by
\[
(f\cdot a)(x)=f(ax)\qquad\text{for all }x\in A.
\]
The first condition needed to construct the Arens multiplication on
$A_{n}^{\sim}$ is
\begin{equation}
f\cdot a\in V\qquad\text{for all }a\in A\text{ and }f\in V. \label{11}%
\end{equation}
When (\ref{11}) holds, for $f\in V$ and $F\in V_{n}^{\sim}$, we may define a
linear functional $F\cdot f$ on $A$ by
\[
(F\cdot f)(a)=F(f\cdot a)\qquad\text{for all }a\in A.
\]
The second condition is%
\begin{equation}
F\cdot f\in V\qquad\text{for all }f\in V\text{ and }F\in V_{n}^{\sim}.
\label{22}%
\end{equation}
When both (\ref{11}) and (\ref{22}) hold, the Arens product of $F,G\in
V_{n}^{\sim}$ is defined by
\[
(F\cdot G)(f)=F(G\cdot f)\qquad\text{for all }f\in V.
\]
Thus the $f$-algebra structure on $V_{n}^{\sim}$ rests on (\ref{11}) and
(\ref{22}).

We study these conditions for $f$-algebras equipped with a Riesz norm. Our
first result (Theorem \ref{thm:first}) shows that (\ref{11}) holds for
$V=A^{\ast}$ if and only if $A$ is bounded, that is, for every $a\in A_{+}$
there exists $\mu>0$ such that $a^{2}\leq\mu a$. Our main positive result
(Theorem \ref{thm:main}) states that if the Riesz norm of $A$ is
submultiplicative, then both (\ref{11}) and (\ref{22}) hold for $V=A^{\ast}$,
so that $(A^{\ast})_{n}^{\sim}$ is an $f$-algebra under the Arens
multiplication; submultiplicativity automatically implies that $A$ is bounded.
When $A$ is norm complete, $A^{\ast}$ coincides with $A^{\sim}$ and the
situation is covered by Theorem 4.4 of \cite{HdP} (Remark \ref{rem:banach}),
so the interest of Theorem \ref{thm:main} lies in non-complete normed $f$-algebras.

Finally, Example \ref{ex:counter} shows that submultiplicativity cannot simply
be dropped: we construct a bounded $f$-algebra of finitely supported real
sequences, equipped with a non-complete Riesz norm that is not
submultiplicative, for which (\ref{22}) fails. Hence, even when (\ref{11})
holds, (\ref{22}) may fail, and the Arens multiplication on $(A^{\ast}%
)_{n}^{\sim}$ need not exist without an additional assumption.

We assume that the reader is familiar with the theory of Riesz spaces as
presented in \cite{AB,M}, to which we refer for terminology and standard results.

\section{Preliminaries and the first invariance condition}

A real associative algebra $A$ is called a \emph{lattice ordered algebra} (or
a \emph{Riesz algebra}) if it is simultaneously a Riesz space whose positive
cone $A_{+}=\{a\in A:0\leq a\}$ is closed under multiplication; equivalently,
$\left\vert ab\right\vert \leq\left\vert a\right\vert \left\vert b\right\vert
$ for all $a,b\in A$. A lattice ordered algebra is called an $f$%
\emph{-algebra} if, for every $a\in A$, the disjoint complement
\[
\{a\}^{d}=\{x\in A:\left\vert a\right\vert \wedge\left\vert x\right\vert =0\}
\]
is a two-sided ring ideal of $A$. We shall use the following standard facts
(see \cite[Section 2.3]{AB} or \cite[Section 3.1]{M}). In an $f$-algebra $A$,
the equalities $|ab|=|a||b|$ and $a^{2}=|a|^{2}$ hold for all $a,b\in A$, and
$A$ is commutative if it is Archimedean.

Throughout, $A$ is an $f$-algebra equipped with a Riesz (or lattice) norm,
that is,
\[
\left\Vert a\right\Vert \leq\left\Vert b\right\Vert \quad\text{whenever
}a,b\in A\text{ and }\left\vert a\right\vert \leq\left\vert b\right\vert .
\]
All norms under consideration are denoted by the same symbol, the intended
meaning being clear from the context. Such an $A$ is Archimedean and hence
commutative. Unless otherwise stated, no compatibility between the
multiplication and the norm is assumed. The norm dual $A^{\ast}$ is a Banach
lattice and an order ideal of $A^{\sim}$ (see \cite[Theorem 3.49]{AB} and
\cite[Proposition 1.3.7]{M}); it separates the points of $A$.

Condition (\ref{11}) always holds for $V=A^{\sim}$, as in \cite{BH,HdP}. To
characterize it for $V=A^{\ast}$ we use the following notion from \cite{BS}.

\begin{definition}
\emph{The }$f$\emph{-algebra }$A$\emph{ is }bounded\emph{ if for every }$a\in
A$\emph{ there exists }$\mu>0$\emph{ such that} $a^{2}\leq\mu\left\vert
a\right\vert $.
\end{definition}

Since $a^{2}=\left\vert a\right\vert ^{2}$, this is equivalent to requiring
that for every $a\in A_{+}$ there is $\mu>0$ with $a^{2}\leq\mu a$.

We start with two lemmas. The first is a form of the fact that a point
separating ideal of $A^{\sim}$ detects the order of $A$.

\begin{lemma}
\label{lem:pos}Let $A$ be an $f$-algebra, and let $V$ be an order ideal of
$A^{\sim}$ that separates the points of $A$. If $a\in A$ and $f(a)\geq0$ for
every $f\in V_{+}$, then $a\in A_{+}$.
\end{lemma}

\begin{proof}
Suppose that $a^{-}\neq0$. Since $V$ separates the points of $A$, there is
$h\in V$ with $h(a^{-})\neq0$, and then $h^{+}(a^{-})>0$ or $h^{-}(a^{-})>0$.
Thus, there exists $f\in V_{+}$ such that $f\left(  a^{-}\right)  >0$. From
\cite[Theorem 1.23]{AB} it follows that there is $g\in A^{\sim}$ with $0\leq
g\leq f$ and $g(-a)=f(a^{-})$. But then $g\in V_{+}$ because $V$ is an order
ideal of $A^{\sim}$. Hence $g\left(  a\right)  \geq0$ by hypothesis and
$f\left(  a^{-}\right)  =-g\left(  a\right)  \leq0$, a contradiction.
\end{proof}

The second is the following estimate, which is implicit in \cite{HdP}.

\begin{lemma}
\label{lem:mod}Let $A$ be an $f$-algebra. If $f\in A^{\sim}$ and $a\in A$,
then%
\[
\left\vert f\cdot a\right\vert (x)\leq\left\vert f\right\vert (\left\vert
a\right\vert x)\quad\text{for all }x\in A_{+}.
\]

\end{lemma}

\begin{proof}
This follows straightforwardly from the Riesz--Kantorovich formula
\cite[Theorem 1.18]{AB}.
\end{proof}

We are in position now to prove the central result of this section.

\begin{theorem}
\label{thm:first}Let $A$ be an $f$-algebra with a Riesz norm. Then
\emph{(\ref{11})} holds for $V=A^{\ast}$ if and only if $A$ is bounded.
\end{theorem}

\begin{proof}
Assume that $f\cdot a\in A^{\ast}$ for all $a\in A$ and $f\in A^{\ast}$. Fix
$a\in A_{+}$ and define the map $T_{a}:A^{\ast}\rightarrow A^{\ast}$ by%
\[
T_{a}(f)=f\cdot a\quad\text{for all }f\in A^{\ast}.
\]
It is an orthomorphism of $A^{\ast}$ (see \cite[p.\ 49]{HdP}). Since $A^{\ast
}$ is a Banach lattice, every orthomorphism of $A^{\ast}$ is central by
Wickstead's theorem \cite[Corollary 4.2]{W} (see also \cite[Theorem 3.29]%
{AA}), so there is $\mu>0$ such that
\[
\left\vert T_{a}(f)\right\vert \leq\mu\left\vert f\right\vert \quad\text{for
all }f\in A^{\ast}.
\]
For $0\leq f\in A^{\ast}$ this gives $f\cdot a\leq\mu f$, and evaluating at
$a\geq0$ we get $f(a^{2})\leq\mu f(a)$, i.e.\ $f(\mu a-a^{2})\geq0$ for all
$f\in\left(  A^{\ast}\right)  _{+}$. By Lemma \ref{lem:pos}, $a^{2}\leq\mu a$
and so $A$ is bounded.

Conversely, let $A$ be bounded, $a\in A$ and $f\in A^{\ast}$. By \cite[Lemma
2.2]{BS} there exists $\mu>0$ with $\left\vert ax\right\vert \leq\mu\left\vert
x\right\vert $ for all $x\in A$. As $\left\vert ax\right\vert =\left\vert
a\right\vert \left\vert x\right\vert $, this means $\left\vert a\right\vert
x\leq\mu x$ for all $x\in A_{+}$. By Lemma \ref{lem:mod}, for $x\in A_{+}$,%
\[
\left\vert f\cdot a\right\vert (x)\leq\left\vert f\right\vert (\left\vert
a\right\vert x)\leq\left\vert f\right\vert (\mu x)=\mu\left\vert f\right\vert
(x).
\]
Thus $\left\vert f\cdot a\right\vert \leq\mu\left\vert f\right\vert $. Since
$A^{\ast}$ is an order ideal of $A^{\sim}$ and $\left\vert f\right\vert \in
A^{\ast}$, we conclude that $f\cdot a\in A^{\ast}$ and the proof is complete.
\end{proof}

\begin{remark}
\label{rem:ideals}The proof of sufficiency does not use the norm. If $A$ is a
bounded Archimedean $f$-algebra, then $\left\vert f\cdot a\right\vert \leq
\mu\left\vert f\right\vert $ for all $f\in A^{\ast}$ and $a\in A$ \emph{(}with
$\mu$ depending on $a$\emph{)}, so \emph{(\ref{11})} holds for every order
ideal $V$ of $A^{\ast}$. Thus, for bounded $f$-algebras, \emph{(\ref{11})} is
never an obstruction, and the only real issue is \emph{(\ref{22})}.
\end{remark}

\section{The main result and the counterexample}

Let $A$ be bounded. By Theorem \ref{thm:first}, $f\cdot a\in A^{\ast}$ for all
$a\in A$ and $f\in A^{\ast}$, so for $F\in(A^{\ast})_{n}^{\sim}$ the linear
functional $F\cdot f$ given by%
\[
(F\cdot f)(a)=F(f\cdot a)\quad\text{for }a\in A
\]
is well defined. As in \cite[proof of Theorem 4.1]{HdP}, $F\cdot f\in A^{\sim
}$. It is natural to ask when $F\cdot f\in A^{\ast}$, that is, when (\ref{22})
holds for $V=A^{\ast}$. The following theorem gives a sufficient condition.

\begin{theorem}
\label{thm:main}Let $A$ be an $f$-algebra with a submultiplicative Riesz norm.
Then $A$ is bounded, conditions \emph{(\ref{11})} and \emph{(\ref{22})} hold
for $V=A^{\ast}$, and $(A^{\ast})_{n}^{\sim}$ is an $f$-algebra with respect
to the Arens multiplication
\[
(F\cdot G)(f)=F(G\cdot f)\quad\text{for all }F,G\in(A^{\ast})_{n}^{\sim}\text{
and }f\in A^{\ast}.
\]

\end{theorem}

\begin{proof}
Let $a,x\in A$ and $f\in A^{\ast}$. Since the norm is submultiplicative,
\[
\left\vert (f\cdot a)(x)\right\vert =\left\vert f(ax)\right\vert
\leq\left\Vert f\right\Vert \left\Vert ax\right\Vert \leq\left\Vert
f\right\Vert \left\Vert a\right\Vert \left\Vert x\right\Vert ,
\]
so $f\cdot a\in A^{\ast}$ and
\[
\left\Vert f\cdot a\right\Vert \leq\left\Vert f\right\Vert \left\Vert
a\right\Vert .
\]
Thus (\ref{11}) holds for $V=A^{\ast}$, and $A$ is bounded by Theorem
\ref{thm:first}.

Now let $F\in(A^{\ast})_{n}^{\sim}$. Being order bounded, $F$ is the
difference of two positive linear functionals on the Banach lattice $A^{\ast}%
$, each of which is norm continuous (see \cite[Theorem 4.3]{AB}); so $F$ is
norm continuous. Therefore, for all $a\in A$,
\[
\left\vert (F\cdot f)(a)\right\vert =\left\vert F(f\cdot a)\right\vert
\leq\left\Vert F\right\Vert \left\Vert f\cdot a\right\Vert \leq\left\Vert
F\right\Vert \left\Vert f\right\Vert \left\Vert a\right\Vert .
\]
Hence $F\cdot f\in A^{\ast}$, which is (\ref{22}).

With (\ref{11}) and (\ref{22}) available, the Arens multiplication on
$(A^{\ast})_{n}^{\sim}$ is well defined, and the proof that it makes
$(A^{\ast})_{n}^{\sim}$ an $f$-algebra follows the proof of \cite[Theorem
4.4]{HdP}; see also \cite[Theorems 3.1 and 3.4]{BH}.
\end{proof}

\begin{remark}
\label{rem:banach}

\begin{enumerate}
\item[\emph{(i)}] In \emph{Theorem \ref{thm:main}}, submultiplicativity can be
replaced by the weaker condition that there exists $C>0$ such that $\left\Vert
ab\right\Vert \leq C\left\Vert a\right\Vert \Vert b\Vert$ for all $a,b\in A$.
Indeed, in that case, $C\left\Vert .\right\Vert $ is a submultiplicative Riesz
norm equivalent to $\left\Vert .\right\Vert $. Equivalent norms give the same
dual $A^{\ast}$, so the theorem applies.

\item[\emph{(ii)}] If $A$ is norm complete, then $A^{\ast}=A^{\sim}$
\emph{(}see, e.g., \emph{\cite[Corollary 4.5]{AB})}. We are then in the case
$V=A^{\sim}$, which is covered by \emph{\cite{HdP}}. Hence any counterexample
to the assertion of \emph{\cite{HdP}} for $V=A^{\ast}$ must involve a
non-complete Riesz norm, as in the example below.
\end{enumerate}
\end{remark}

Without submultiplicativity of the norm, (\ref{22}) may fail, as the following
example shows.

\begin{example}
\label{ex:counter}Let $A=c_{00}$ be the bounded $f$-algebra of all finitely
supported real sequences \emph{(}indexed by $\mathbb{N}=\left\{
1,2,....\right\}  $\emph{)}. For each $k\in\mathbb{N}$, define $e_{k}\in A$ by%
\[
e_{k}\left(  n\right)  =\delta_{k,n}\text{\quad for all }n\in\mathbb{N}.
\]
Equip $A$ with the Riesz norm defined by%
\[
\left\Vert a\right\Vert =\sum_{n=1}^{\infty}\frac{\left\vert a\left(
n\right)  \right\vert }{2^{n}}\quad\text{for all }a\in A.
\]
Observe that if $f\in A^{\ast}$ and $n\in\mathbb{N}$ then%
\[
\left\vert 2^{n}f\left(  e_{n}\right)  \right\vert \leq2^{n}\left\Vert
f\right\Vert \left\Vert e_{n}\right\Vert =\left\Vert f\right\Vert .
\]
It follows that the map $F:A^{\ast}\rightarrow\mathbb{R}$ given by%
\[
F\left(  f\right)  =\sum_{n=1}^{\infty}\frac{2^{n}f\left(  e_{n}\right)
}{n^{2}}\quad\text{for all }f\in A^{\ast}%
\]
is well defined. Actually, a short moment's thought reveals that $0\leq
F\in\left(  A^{\ast}\right)  _{n}^{\sim}$.

Now, let $f\in A^{\ast}$ be defined by%
\[
f\left(  e_{n}\right)  =\frac{1}{2^{n}}\quad\text{for all }n\in\mathbb{N}.
\]
For every $n\in\mathbb{N}$, we have%
\[
\left(  F\cdot f\right)  \left(  e_{n}\right)  =F\left(  f\cdot e_{n}\right)
=\sum_{k=1}^{\infty}\frac{2^{k}\left(  f\cdot e_{n}\right)  \left(
e_{k}\right)  }{k^{2}}=\sum_{k=1}^{\infty}\frac{2^{k}f\left(  e_{n}%
e_{k}\right)  }{k^{2}}=\frac{1}{n^{2}}.
\]
Consequently,%
\[
\left\Vert 2^{n}e_{n}\right\Vert =1\quad\text{and\quad}\left(  F\cdot
f\right)  \left(  2^{n}e_{n}\right)  =\frac{2^{n}}{n^{2}},
\]
showing that $F\cdot f\notin A^{\ast}$. Thus \emph{(\ref{22})} fails for
$V=A^{\ast}$, so $\left(  A^{\ast}\right)  _{n}^{\sim}$ cannot be equipped
with the Arens multiplication.
\end{example}

\begin{remark}

\begin{enumerate}
\item[\emph{(i)}] The norm of \emph{Example \ref{ex:counter}} is not even
$C$-submultiplicative: $\left\Vert e_{n}\right\Vert =2^{-n}$, whereas
$\left\Vert e_{n}\right\Vert ^{2}=4^{-n}$. This is consistent with
\emph{Remark \ref{rem:banach} (i)}. Moreover, the completion of $A$ is the
space of sequences $a$ with $\sum_{n}2^{-n}|a(n)|<\infty$, which is not closed
under coordinatewise multiplication \emph{(}consider $a(n)=2^{n/2}$\emph{)};
so the multiplication of $A$ does not extend to the completion, and the
situation of \emph{Remark \ref{rem:banach} (iii)} cannot be recovered by completing.

\item[\emph{(ii)}] Let $B$ be the $f$-algebra of all eventually constant real
sequences with the same norm. It is a bounded, non-complete, Riesz-normed
$f$-algebra with identity in which $A$ is norm dense. Therefore, $B^{\ast}$ is
canonically the same space $A^{\ast}$, and the same $F$ and $f$ satisfy
$F\cdot f\notin B^{\ast}$. So the failure of \emph{(\ref{22})} is not related
to the absence of an identity.

\item[\emph{(iii)}] Conceptually, when \emph{(\ref{11})} holds,
\emph{(\ref{22})} is a continuity property of the maps $\Phi_{f}:A\rightarrow
A^{\ast}$ defined by $\Phi_{f}(a)=f\cdot a$ for all $a\in A$. If the norm is
submultiplicative, then $\left\Vert \Phi\right\Vert \leq\left\Vert
f\right\Vert $. In the example, $\left\Vert \Phi_{f}\left(  e_{n}\right)
\right\Vert =\left\Vert 2^{-n}e_{n}\right\Vert =1$ while $\left\Vert
e_{n}\right\Vert =2^{-n}$, so $\Phi_{f}$ is unbounded, although every
multiplication operator by an element is norm continuous \emph{(}as $A$ is
bounded\emph{)}. Thus separate continuity of the multiplication is not enough,
while the joint continuity of \emph{Remark \ref{rem:banach} (i)} is. We do not
know whether joint continuity is also necessary for \emph{(\ref{22})}.
\end{enumerate}
\end{remark}

\section*{Acknowledgement}

The authors thank Professor Ben de Pagter for his kind response to their
inquiry regarding the remark following Theorem 4.4 of \cite{HdP}. Any
inaccuracies in this note are the authors' own.

\end{document}